\documentclass[12pt,a4paper]{amsart}
\usepackage{amssymb}
\usepackage{mathrsfs}
\headsep=1cm \numberwithin{equation}{section}
\newtheorem{thm}{Theorem}[section]
\newtheorem{cor}[thm]{Corollary}
\newtheorem{lem}[thm]{Lemma}
\newtheorem{prop}[thm]{Proposition}
\theoremstyle{definition}

\theoremstyle{remark}

\numberwithin{equation}{section}

\newcommand\Supp{\operatorname{Supp}}
\newcommand\Ass{\operatorname{Ass}}
\newcommand\mAss{\operatorname{mAss}}

\newcommand\Spec{\operatorname{Spec}}
\newcommand\Rad{\operatorname{Rad}}
\newcommand\Hom{\operatorname{Hom}}
\newcommand\Ext{\operatorname{Ext}}

\newcommand\height{\operatorname{height}}

\newcommand\Max{\operatorname{Max}}
\begin{document}\title[Bass numbers and Cominimaxness]
{On the finiteness of Bass numbers of local cohomology modules and Cominimaxness}
\author[Kamal Bahmanpour,  Reza Naghipour and Monireh Sedghi]{Kamal Bahmanpour,  Reza Naghipour$^*$ and Monireh Sedghi\\\\\\\,
\vspace*{0.5cm}Dedicated to Professor Hyman Bass}
\vspace*{0.5cm}
\address{Department of Mathematics, Ardabil Branch, Islamic Azad University, Ardabil, Iran.}
\email{bahmanpour.k@gmail.com}
\address{Department of Mathematics, University of Tabriz, Tabriz, Iran;
and School of Mathematics, Institute for Research in Fundamental
Sciences (IPM), P.O. Box: 19395-5746, Tehran, Iran.}
\email{naghipour@ipm.ir} \email {naghipour@tabrizu.ac.ir}

\address{Department of Mathematics, Azarbaijan University of Shahid Madani, Tabriz, Iran. }%
\email{sedghi@azaruniv.edu}%

\thanks{ 
This research was in part supported by a grant from IPM (No. 89130053, 89130048).\\
$^*$Corresponding author: e-mail: {\it naghipour@ipm.ir} (Reza Naghipour)}%
\begin{abstract}
In this paper, we continue the study of cominimaxness modules with
respect to an ideal of a commutative Noetherian ring (cf.  \cite{ANV}), and
 Bass numbers of local cohomology modules.
  Let $R$ denote a commutative Noetherian local ring and $I$ an ideal
of $R$. We first show that the Bass numbers $\mu^0(\frak p, H^2_I(R))$ and
$\mu^1(\frak p, H^2_I(R))$ are finite for all $\frak p\in \Spec R$,
whenever $R$ is regular. As a consequence, it follows that the
Goldie dimension of $H^2_I(R)$ is finite. Also, for a finitely
generated $R$-module $M$ of dimension $d$, it is shown that the Bass
numbers of $H^{d-1}_{I}(M)$ are finite if and only if $\Ext^i_R(R/I,
H^{d-1}_{I}(M))$  be minimax for all $i\geq0$. Finally, we prove that
if $\dim R/I=2$, then the Bass numbers of $H^{n}_{I}(M)$ are finite
if and only if $\Ext^i_R(R/I, H^{n}_{I}(M))$ be minimax, for all
$i\geq0$, where $n$ is a non-negative integer.
\end{abstract}
\maketitle
\section{Introduction}
Throughout this paper, let $R$ denote a commutative Noetherian ring
(with identity) and $I$ an ideal of $R$. For an $R$-module $M$, the
$i^{th}$ local cohomology module of $M$ with respect to $I$ is
defined as$$H^i_I(M) = \underset{n\geq1} {\varinjlim}\,\,
\text{Ext}^i_R(R/I^n, M).$$ We refer the reader to \cite{Gr1} and
\cite{BS} for more details about local cohomology. An important
problem in commutative algebra is to determine when the Bass numbers
of the $i^{th}$ local cohomology module $H^i_I(M)$ is finite. In
\cite{Hu} Huneke conjectured that for any ideal $I$ in a regular
local ring $(R, \frak m, k)$, the Bass numbers $$\mu^i(\frak p,
H^i_I(R))= \dim _{k(\frak p)}\Ext^i_{R_{\frak p}}(k(\frak p),
H^i_{IR_{\frak p}}(R_{\frak p}))$$ are finite for all $i$ and $j$,
where $k(\frak p):=R_{\frak p}/\frak p R_{\frak p}$. In particular
the injective resolution of the local cohomology has only finitely
many copies of the injective hull of $R/\frak p$ for any $\frak p$.

There is evidence that this conjecture is true. It is shown that by
Huneke and Sharp  \cite{HS} and Lyubeznik \cite{Ly1, Ly2} that the
conjecture holds for a regular local ring containing  a field. We
remark that the Bass numbers might be infinite if $R$ is not
regular. For example, if $I:=(x, y)R\subseteq R:=k[x, y, z,
w]/(xz-yw)$, then $\mu^0(\frak m, H^2_I(R))=\infty$ for $\frak
m:=(x, y, z, w)R $ (see \cite{Ha}).

 We say that $M$ is a minimax module if there is a
finitely generated submodule $N$ of $M$, such that $M/N$ is
Artinian. The interesting class of minimax modules was introduced by
H. Z\"{o}shinger in \cite{Zo1} and he has in \cite{Zo1} and
\cite{Zo2} given many equivalent conditions for a module to be
minimax. Finally, the $R$-module $M$ is said to be an $I$-cominimax
if support of $M$ is contained in $V(I)$ and $\Ext^i_R(R/I, M)$ is
minimax for all $i\geq0$. The concept of the $I$-cominimax modules
were introduced in \cite{ANV} as a generalization of important
notion of $I$-cofinite modules.

The main subject of the paper is to continue the study of $I$-cominimaxness
properties and the finiteness of Bass numbers of local cohomology
modules. First, we  provide a partial answer to
Huneke's conjecture. Namely, it will be shown that the Bass numbers
of $\mu^0(\frak p, H^2_I(R))$ and  $\mu^1(\frak p, H^2_I(R))$ are
finite for all ideals $I$ in a regular local ring $R$ and for all
$\frak p\in \Spec R$. Pursuing this point of view further we
establish the $R$-module $H^2_I(R)$ has finite Goldie dimension.
Recall that an $R$-module $M$ is said to have finite Goldie
dimension if $M$ does not contain an infinite direct sum of non-zero
submodules, or equivalently the injective hull $E(M)$ of $M$
decompose as a finite direct sum of indecomposable (injective)
submodules.

Since the Bass numbers of the local cohomology modules $H^i_I(M)$ are not finite in general,
several attempts made to find some conditions for the
ideal $I$ to have finiteness for the Bass numbers of the local
cohomology modules with support in $I$ in terms of $I$-cominimaxness
properties. The following result is in this way:

\begin{thm}
Let $(R,\mathfrak{m})$ be a local (Noetherian) ring, $I$ an ideal of
$R$ and $M$ a finitely generated $R$-module of dimension d. Then the
Bass numbers of $H^{d-1}_{I}(M)/\Gamma _{\frak m}(H^{d-1}_{I}(M))$
are finite.
\end{thm}

One of our tools for proving Theorem 1.1 is the following:

\begin{prop}
Let $(R,\mathfrak{m})$, $I$ and $M$ be as in Theorem 1.1. Then
$H^1_{\frak m}(H^{d-1}_{I}(M))$ is an Artinian $R$-module, where
$d=\dim M$.
\end{prop}

Pursuing this point of view further we derive the following
consequence of Theorem 1.1.

\begin{thm}
Let $(R,\mathfrak{m})$ be a local (Noetherian) ring, $I$ an ideal of
$R$ and $M$ a finitely generated $R$-module of dimension $d$. Then
the following conditions are equivalent:

{\rm(i)} The Bass numbers of $H^{d-1}_{I}(M)$ are finite.

{\rm(ii)} ${\rm Soc}(H^{d-1}_{I}(M))$ is finitely generated.

{\rm (iii)} $H^{d-1}_{I}(M)$ is $I$-cominimax.
\end{thm}

As the second main result of this paper we obtain a characterization of the finiteness
of Bass numbers of $i^{th}$ local cohomology modules of $M$ with respect to $I$ of
dimension 2, i.e. $\dim R/I=2$.
More precisely we shall show that:

\begin{thm}
Let $(R,\mathfrak{m})$ be a local (Noetherian)  ring and let $I$ be
an ideal of $R$ with $\dim R/I=2$. Let $M$ be a finitely generated
$R$-module and $i$ a non-negative integer. Then the Bass numbers of
$H^{i}_{I}(M)$ are finite if and only if the $R$-module
$H^{i}_{I}(M)$ is $I$-cominimax.
\end{thm}

The proof of Theorem 1.4 is given in 2.17. As an application, we
derive the following consequences of Theorem 1.4, which is a
characterization of $I$-cominimaxness of $H^{i}_{I}(M)$ in terms of
finiteness of Bass numbers of local cohomology modules
$H^{i}_{I}(M)$ for certain ideal $I$ of $R$.

\begin{thm}
Let $(R,\mathfrak{m})$ be a local (Noetherian) ring, and let $I$ be
an ideal of $R$ such that $\dim R/I=2$. Let $M$ be a finitely
generated $R$-module and $t$ a non-negative integer. Then the
following conditions are equivalent:

{\rm(i)} ${\rm Soc}\, H^{i}_{I}(M)$ is finitely generated for all
$i\leq t$.

{\rm(ii)} $\Ext^{j}_{R}(R/\mathfrak{m},H^{i}_{I}(M))$ is finitely
generated for all $j \geq 0$ and all $i\leq t-1$.

{\rm(iii)} $H^{i}_{I}(M)$ is $I$-cominimax for all $i\leq t-1$.
\end{thm}

Using Theorem 1.5 we obtain some results as following:
\begin{cor}
Let $(R,\mathfrak{m})$ be a regular local ring of dimension $d\geq
2$ containing a field and let $I$ be an ideal of $R$ with $\dim
R/I=2$.
 Then the $R$-module $H^{i}_{I}(R)$ is $I$-cominimax for all $i\geq0$.
\end{cor}

\begin{cor}
Let $R$ be a Noetherian ring, $M$ a finitely generated $R$-module
and $I$  an ideal of $R$ such that $\dim M/IM=2$. Let
$$\Sigma:=\{\frak p \in \Spec R:\,\,\mu^j(\frak p,H^{i}_{I}(M)
)=\infty\,\,{\rm for}\,\,{\rm some}\,\,{\rm integers}\,\,i\geq
0\,\,{\rm and}\,\,j\geq 0 \}.$$ Then $\Sigma$ is  countable  and
$\Sigma \subseteq {\rm Max}(R).$
\end{cor}
\begin{cor}
Let $(R,\mathfrak{m})$ be a local (Noetherian) ring, $M$  a
finitely generated $R$-module and $I$ an ideal of $R$  such that
$\dim M/IM=3$. Then the set
$$\Sigma:=\{\frak p \in \Spec R:\,\,\mu^j(\frak p,H^{i}_{I}(M)
)=\infty\,\,{\rm for}\,\,{\rm some}\,\,{\rm integers}\,\,i\geq
0\,\,{\rm and}\,\,j\geq 0 \},$$ is  countable.
\end{cor}

Throughout this paper, $R$ will always be a commutative Noetherian
ring with non-zero identity and $I$ will be an ideal of $R$.
For each $R$-module $L$, we denote by
 ${\rm Ass h}_RL$ (resp. ${\rm mAss}_RL$) the set $\{\frak p\in \Ass
_RL:\, \dim R/\frak p= \dim L\}$ (resp. the set of minimal primes of
$\Ass_RL$). We shall use $\Max R$ to denote the set of all maximal
ideals of $R$. Also, for any ideal $\frak a$ of $R$, we denote
$\{\frak p \in {\rm Spec}\,R:\, \frak p\supseteq \frak a \}$ by
$V(\frak a)$. Finally, for any ideal $\frak{b}$ of $R$, {\it the
radical of} $\frak{b}$, denoted by $\Rad(\frak{b})$, is defined to
be the set $\{x\in R \,: \, x^n \in \frak{b}$ for some $n \in
\mathbb{N}\}$. For any unexplained notation and terminology we refer
the reader to \cite{BS} and \cite{Mat}.

\section{The Results}

Let $(R,\mathfrak{m})$ be a local (Noetherian) ring and $I$ an ideal of $R$.
As noted in the introduction, it is well known that the Bass numbers $\mu^i(\frak p, H^2_I(R))$
are finite whenever $R$ is a regular ring containing a field. Our first result shows that
if $R$ is regular, then $\mu^i(\frak p, H^2_I(R))$ is finite for $i\in \{0, 1\}$.

\begin{thm}
Let $(R,\mathfrak{m})$ be a regular local ring and $I$ an ideal of
$R$. Then the Bass numbers $\mu^{0}(\mathfrak{p},H^{2}_{I}(R))$ and
 $\mu^{1}(\mathfrak{p},H^{2}_{I}(R))$ are finite for all
  $\mathfrak{p}\in \Spec R.$
\end{thm}
\proof In view of  \cite[Proposition 2.20]{BN3} the
$R_{\mathfrak{p}}$-modules
$\Ext^{j}_{R_{\mathfrak{p}}}(\kappa(\mathfrak{p}),H^{i}_{IR_{\mathfrak{p}}}(R_{\mathfrak{p}}))$
are finitely generated for all $j\geq0$ and $i=0,1$. Hence it
follows from \cite[Corollary 3.5]{Kh} that
$\mu^{0}(\mathfrak{p},H^{2}_{I}(R))$ and
$\mu^{1}(\mathfrak{p},H^{2}_{I}(R))$ are finite.
\qed\\

As a corollary of the previous theorem, we deduce that the Goldie dimension of $H^{2}_{I}(R)$
is finite.

\begin{cor}
Let $(R,\mathfrak{m})$ be a regular local ring and $I$ an ideal of
$R$. Then the $R$-module $H^{2}_{I}(R)$ has finite Goldie dimension.
\end{cor}

\proof As $Ass_{R}(H^{2}_{I}(R))$ is finite, by \cite[Theorem
2.4]{BN2}, it follows from Theorem 2.1 and the definition that
$H^{2}_{I}(R)$ has finite Goldie dimension.\qed\\

For an arbitrary $R$-module $M$, the next reult gives us a necessary and
sufficient condition for the finiteness of $\mu^{i}(\mathfrak{p}, M)$, the Bass
numbers of $M$, in terms of the minimaxness of $\Ext^{i}_{R}(R/I,M)$. Recall
that for any module $L$ over a local ring $(R, \frak m)$, the socle of $L$, denoted
by $Soc(L)$, is defined to be the $R$-module $\Hom_R(R/\frak m, L)$.
\begin{thm}
Let $(R,\mathfrak{m})$ be a local (Noetherian) ring and $I$ an ideal
of $R$ such that $\dim R/I=1$. Suppose that $M$ is an $R$-module and
$n\geq 0$ an integer. Then the following conditions are equivalent:

{\rm(i)} $\mu ^{i}(\mathfrak{p},M)$ is finite for all
$\mathfrak{p}\in V(I)$ and for all $i\leq n$.

{\rm(ii)} $\Ext^{i}_{R}(R/I,M)$ is minimax for all $i\leq n$.
\end{thm}

\proof First we show (i)$\Longrightarrow$(ii). Let $\mathfrak{q}\in
\mAss_{R}(R/I)$ and $i\leq n$. Then $\mathfrak{q}\neq \mathfrak{m}$,
and so $\Rad(\mathfrak{q}+Rx)=\mathfrak{m}$ for all $x\in
\mathfrak{m}\backslash \mathfrak{q}$. Hence it follows from
\cite[Proposition 1]{DM} that the $R$-module
$\Ext^{i}_{R}(R/\mathfrak{q}+Rx,M)$ is finitely generated for all
$i\leq n$. (Note that $\Ext^{i}_{R}(R/\mathfrak{m},M)$ is finitely
generated for all $i\leq n$.) Now, for $x\in \mathfrak{m}\backslash
\mathfrak{q}$ the exact sequence$$0\rightarrow
R/\mathfrak{q}\stackrel{x} \longrightarrow
R/\mathfrak{q}\longrightarrow R/\mathfrak{q}+Rx\rightarrow 0,$$
provides the following exact
sequence:$$\Ext^{i}_{R}(R/\mathfrak{q}+Rx,M)  \longrightarrow
\Ext^{i}_{R}(R/\mathfrak{q},M)\stackrel{x} \longrightarrow
\Ext^{i}_{R}(R/\mathfrak{q},M),$$ which implies that the $R$-module
${\rm Soc}(\Ext^{i}_{R}(R/\mathfrak{q},M))$ is finitely generated,
 for all $i\leq n$. On the other hand, if $L:=\Ext^{i}_{R}(R/\mathfrak{q},M)$, then as by assumption
 the $R_{\mathfrak{q}}$-module $L_{\mathfrak{q}}$ is finitely generated, it follows that there exists
 a finitely generated submodule $K$ of $L$, such that $(L/K)_{\mathfrak{q}}=0$, and so
 $\Supp (L/K) \subseteq V(\mathfrak{m})$. Finally, the exact sequence
 $$0\rightarrow K \longrightarrow L \longrightarrow L/K \longrightarrow 0.$$induces the
 exact sequence $${\rm Soc}(L)\longrightarrow {\rm Soc}(L/K) \longrightarrow \Ext^{1}_{R}(R/\mathfrak{m},K).$$
  Hence the $R$-module ${\rm Soc}(L/K)$ is finitely generated.
  Consequently in view of \cite[Proposition 4.1]{Me3} the $R$-module
  $L/K$ is Artinian, that is $L$ is a minimax $R$-module. Now, the assertion
  follows from \cite[Corollary 2.8]{ANV}. \\
In order to prove (ii)$\Longrightarrow$(i), use \cite[Corollary
2.8]{ANV} and the
definition of the minimax modules.\qed\\

\begin{cor}
Let $R$, $I$ and $M$ be as in Theorem {\rm 2.3.} Then the following
conditions are equivalent:

{\rm (i)} $\mu ^{i}(\mathfrak{p},M)$ is finite for all $i\geq0$ and
all primes $\mathfrak{p}\in V(I)$.

{\rm (ii)} $\Ext^{i}_{R}(R/I,M)$ is minimax for all $i\geq 0$.

\proof The result follows from Theorem 2.3. \qed\\
\end{cor}

\begin{cor}
Let $(R,\mathfrak{m})$ be a local (Noetherian) ring.  Let $I$ and $J$ be tow
ideals of $R$ such that $\dim R/I\leq 1$ and $\dim R/J=1$. Then the
$R$-modules $\Ext^{i}_{R}(R/J,H^{j}_{I}(R))$ are minimax for all
integers $i,j\geq0$.
\end{cor}

\proof The assertion follows from \cite[Corollary 2.10]{BN3} and Theorem 2.3.\qed\\

Before we state the next corollary, recall that for an ideal $I$ of  a commutative ring
$R$, the ideal transform $D_I(R)$ of $I$ defined by $$D_I(R):=\underset{n\geq1} {\varinjlim}\,\,
\text{Hom}_R(I^n, R).$$  See  \cite[Section 2.2]{BS} for the basic properties of ideal transforms.\\

\begin{cor}
Let $(R,\mathfrak{m})$ be a regular local ring and let $I, J$ be two
ideals of $R$ such that $\dim R/I=1$. Then

{\rm(i)} For any integer $i\geq 0$, the $R$-modules
$\Ext^{i}_{R}(R/I,D_{J}(R))$ and $\Ext^{i}_{R}(R/I,H_{J}^{1}(R))$
are minimax.

{\rm(ii)} The $R$-modules $\Hom_{R}(R/I,H_{J}^{2}(R))$ and
$\Ext_{R}^{1}(R/I,H_{J}^{2}(R))$ are minimax.

{\rm(iii)}  For all integers $i, j\geq 0$, the $R$-module
$\Ext_{R}^{i}(R/I,H_{J}^{j}(R))$ is minimax, whenever $R$ contains a
field.
\end{cor}

\proof (i) follows from \cite[Proposition 2.20]{BN3} and Theorem
2.3.  In order to prove (ii) use Theorems 2.3 and 2.1. Finally, (iii)
follows from Corollary 2.4 and \cite{HS, Ly2}. \qed\\

In order to state next results, we need here some preliminary results about
the Artinianess of local cohomology modules.\\

\begin{thm}
Let $(R,\mathfrak{m})$ be a local (Noetherian) ring and $I$ an ideal
of $R$ such that $\dim R/I=2$. Then for any  finitely generated
$R$-module $M$, the $R$-module $H^{1}_{\mathfrak{m}}(H^{i}_{I}(M))$
is Artinian for all $i\geq 0$.
\end{thm}

\proof
Since $\dim R/I=2$, there exists $x \in \mathfrak{m}$ such that
$\dim R/I+Rx=1$. Let $J:=I+Rx$. As $H^{1}_{\mathfrak{m}}(H^{0}_{I}(M))$ is Artinian,
we may assume that $i\geq 1$.
Then by \cite[Corollary 1.4]{Sc}, there exists the following exact
sequence, $$0\rightarrow H^{1}_{J}(H^{i-1}_{I}(M)) \rightarrow
H^{i}_{J}(M) \rightarrow H^{0}_{J}(H^{i}_{I}(M))\rightarrow 0.$$
Now, as $\dim H^{1}_{J}(H^{i-1}_{I}(M))\leq 1$, it follows from this
exact sequence that the sequence$$H^{1}_{\mathfrak{m}}(H^{i}_{J}(M))
\rightarrow H^{1}_{\mathfrak{m}}(H^{0}_{J}(H^{i}_{I}(M)))\rightarrow
0,$$is exact, and so by \cite[Corollary 2.16]{BN3} the $R$-module $
H^{1}_{\mathfrak{m}}(H^{0}_{J}(H^{i}_{I}(M)))$ is Artinian. On the
other hand, since $\dim R/J=1$ there exists $y \in \mathfrak{m}$
such that $J+Ry$ is $\mathfrak{m}$-primary and so
$$H^{1}_{J+Ry}(H^{0}_{J}(H^{i}_{I}(M)))=H^{1}_{\mathfrak{m}}(H^{0}_{J}(H^{i}_{I}(M))),$$is
Artinian $R$-module. Moreover, using again \cite[Corollary 1.4]{Sc}
it follows that the sequences:$$0\rightarrow H^{1}_{J}(H^{i}_{I}(M))
\rightarrow H^{i+1}_{J}(M) \rightarrow
H^{0}_{J}(H^{i+1}_{I}(M))\rightarrow 0,$$ $$0\rightarrow
H^{1}_{J+Ry}(H^{0}_{J}(H^{i}_{I}(M))) \rightarrow
H^{1}_{\mathfrak{m}}(H^{i}_{I}(M)) \rightarrow
H^{0}_{J+Ry}(H^{1}_{J}(H^{i}_{I}(M)))\rightarrow 0,$$ are exact.
Now, since by \cite[Corollary 2.16]{BN3},
$H^{0}_{J+Ry}(H^{i+1}_{J}(M))$ is Artinian, it follows that
$H^{0}_{J+Ry}(H^{1}_{J}(H^{i}_{I}(M)))$ is Artinian. Consequently
the $R$-module $H^{1}_{\mathfrak{m}}(H^{i}_{I}(M))$ is
Artinian, as required.\qed\\

\begin{thm}
Let $(R,\mathfrak{m})$ be a local (Noetherian) ring, $I$ an ideal of
$R$ and $M$ a finitely generated $R$-module of dimension $d\geq1$.
Then $H^{1}_{\mathfrak{m}}(H^{d-1}_{I}(M))$ is an Artinian
$R$-module.
\end{thm}

\proof We use induction on $d$. When $d=1$, there is nothing to
prove. Now suppose that $d>1$ and the result has been proved for
non-zero finitely generated $R$-modules of dimension $d-1$. By
replacing $M$ by $M/\Gamma_{I}(M)$ we may assume that $M$ is a
non-zero finitely generated $I$-torsion-free $R$-module. Then, by
\cite [Lemma 2.1.1]{BS}, the ideal $I$ contains an element $x$ which
is a non-zerodevisor on $M$. Hence the exact sequence
$$0\longrightarrow M  \stackrel{x} \longrightarrow M \longrightarrow M/xM
\longrightarrow 0$$ induces an exact sequence $$\cdots
\longrightarrow H^{j}_{I}(M)  \stackrel{x} \longrightarrow
H^{j}_{I}(M) \longrightarrow H^{j}_{I}(M/xM) \longrightarrow
H^{j+1}_{I}(M) \stackrel{x} \longrightarrow H^{j+1}_{I}(M)
\longrightarrow \cdots.$$

Therefore we have the following exact sequence, $$0 \longrightarrow
H^{j}_{I}(M)/xH^{j}_{I}(M) \longrightarrow H^{j}_{I}(M/xM)
\longrightarrow (0: _{H^{j+1}_{I}(M)} x) \longrightarrow 0,$$ and so
it follows from $\dim M/xM= d-1$ that
$H^{d-1}_{I}(M)/xH^{d-1}_{I}(M)$ is Artinian. Moreover, by the
inductive hypothesis, the $R$-module
$H^{1}_{\mathfrak{m}}(H^{d-2}_{I}(M/xM))$ is Artinian. Next, in view
of \cite [Corollary 3.3]{NS}, $\Supp H^{d-2}_{I}(M/xM)$ is finite.
Hence we have $\dim H^{d-2}_{I}(M/xM)\leq 1$. Therefore the exact
sequence$$0 \longrightarrow H^{d-2}_{I}(M)/xH^{d-2}_{I}(M)
\longrightarrow H^{d-2}_{I}(M/xM) \longrightarrow (0:
_{H^{d-1}_{I}(M)} x) \longrightarrow 0,$$  provides the exact
sequence
$$H^{1}_{\mathfrak{m}}(H^{d-2}_{I}(M)/xH^{d-2}_{I}(M))
\longrightarrow H^{1}_{\mathfrak{m}}(H^{d-2}_{I}(M/xM))
\longrightarrow H^{1}_{\mathfrak{m}}((0: _{H^{d-1}_{I}(M)} x))
\longrightarrow 0.$$  Hence the $R$-module $H^{1}_{\mathfrak{m}}((0:
_{H^{d-1}_{I}(M)} x))$ is Artinian. Now, from the exact sequences
$$0 \rightarrow (0:_{H^{d-1}_{I}(M)}x) \rightarrow H^{d-1}_{I}(M)  \stackrel{f} \rightarrow xH^{d-1}_{I}(M) \rightarrow 0,$$
$$0 \rightarrow xH^{d-1}_{I}(M)  \stackrel{g}\rightarrow H^{d-1}_{I}(M) \rightarrow H^{d-1}_{I}(M)/xH^{d-1}_{I}(M) \rightarrow 0,$$
we obtain the following exact sequences
$$H^{1}_{\mathfrak{m}}((0:_{H^{d-1}_{I}(M)}x))
\longrightarrow
H^{1}_{\mathfrak{m}}(H^{d-1}_{I}(M))\stackrel{H^{1}_{\mathfrak{m}}(f)}
\longrightarrow H^{1}_{\mathfrak{m}}(xH^{d-1}_{I}(M)),$$
$$H^{0}_{\mathfrak{m}}(H^{d-1}_{I}(M)/xH^{d-1}_{I}(M))
\longrightarrow H^{1}_{\mathfrak{m}}(xH^{d-1}_{I}(M))
\stackrel{H^{1}_{\mathfrak{m}}(g)}\longrightarrow
H^{1}_{\mathfrak{m}}(H^{d-1}_{I}(M)).$$ Then ${\rm
Ker}(H^{1}_{\mathfrak{m}}(f))$ and ${\rm
Ker}(H^{1}_{\mathfrak{m}}(g))$
 are Artinian (note that $H^{1}_{\mathfrak{m}}((0:_{H^{d-1}_{I}(M)}x))$ and $H^{d-1}_{I}(M)/xH^{d-1}_{I}(M)$
 are Artinian), and the sequence $$0 \rightarrow \ker(H^{1}_{\mathfrak{m}}(f)) \rightarrow
 \ker(H^{1}_{\mathfrak{m}}(g)\circ H^{1}_{\mathfrak{m}}(f)) \rightarrow \ker(H^{1}_{\mathfrak{m}}(g)),$$
  is exact. Since $\ker(H^{1}_{\mathfrak{m}}(g\circ f))=(0:_{H^{1}_{\mathfrak{m}}(H^{d-1}_{I}(M))}x)$ and
   $\ker(H^{1}_{\mathfrak{m}}(g)\circ H^{1}_{\mathfrak{m}}(f))$ is Artinian, it follows
   that $(0:_{H^{1}_{\mathfrak{m}}(H^{d-1}_{I}(M))}x)$ is also Artinian.
   Whence according to Melkersson's Theorem \cite[Proposition 1.4]{Me3}
   the $R$-module $H^{1}_{\mathfrak{m}}(H^{d-1}_{I}(M))$ is Artinian,
   and this completes the inductive step.\qed\\

\begin{thm}
Let $(R,\mathfrak{m})$ be a local (Noetherian) ring, $I$ an ideal of
$R$ and $M$ a finitely generated $R$-module of dimension $d$. Then
the Bass numbers of
$H^{d-1}_{I}(M)/\Gamma_{\mathfrak{m}}(H^{d-1}_{I}(M))$ are finite.
\end{thm}

\proof In view of Theorem 2.8, the $R$-module
$H^{1}_{\mathfrak{m}}(H^{d-1}_{I}(M)/\Gamma_{\mathfrak{m}}(H^{d-1}_{I}(M)))$
is Artinian. On the other hand, since $\Supp H^{d-1}_{I}(M)$ is
finite (see \cite [Corollary 3.3]{NS}), it follows that $\dim
H^{d-1}_{I}(M)\leq 1 $, and so for $i\geq 2 $,
$H^{i}_{\mathfrak{m}}(H^{d-1}_{I}(M)/\Gamma_{\mathfrak{m}}(H^{d-1}_{I}(M)))=0$.
Therefore for all $i\neq1$ we have
$H^{i}_{\mathfrak{m}}(H^{d-1}_{I}(M)/\Gamma_{\mathfrak{m}}(H^{d-1}_{I}(M)))=0$.
Consequently, in view of  \cite [Corollary 3.10]{Me3} the
$R$-modules
$\Ext^{i}_{R}(R/\mathfrak{m},H^{d-1}_{I}(M)/\Gamma_{\mathfrak{m}}(H^{d-1}_{I}(M)))$
are finitely generated. (Note that a module is
$\mathfrak{m}$-cofinite if and only if it is Artinian.) Furthermore,
because of $\dim H^{d-1}_{I}(M)\leq 1 $, it follows that $\dim
R/\mathfrak{p}=1$ for all $\mathfrak{p}\in \Supp
H^{d-1}_{I}(M)/\Gamma_{\mathfrak{m}}(H^{d-1}_{I}(M))$ with
$\mathfrak{p}\neq \mathfrak{m}$. Since
$$(H^{d-1}_{I}(M)/\Gamma_{\mathfrak{m}}(H^{d-1}_{I}(M)))_{\mathfrak{p}}\cong
(H^{d-1}_{I}(M))_{\mathfrak{p}},$$for all $\mathfrak{p}\in \Supp
H^{d-1}_{I}(M)/\Gamma_{\mathfrak{m}}(H^{d-1}_{I}(M))$ with
$\mathfrak{p}\neq \mathfrak{m}$, it follows that
$H^{d-1}_{IR_{\mathfrak{p}}}(M_{\mathfrak{p}})\neq0$, and so $\dim
M_{\mathfrak{p}}=d-1$. Consequently
$H^{d-1}_{IR_{\mathfrak{p}}}(M_{\mathfrak{p}})$ is an Artinian
$R_{\mathfrak{p}}$-module. Hence for each $i\geq0$ and each
$\mathfrak{p}\in {\rm Spec} R$, the $i^{th}$ Bass
number$$\mu^{i}(\mathfrak{p},H^{d-1}_{I}(M)/\Gamma_{\mathfrak{m}}(H^{d-1}_{I}(M)))
=\dim_{\kappa(\mathfrak{p})}\Ext^{i}_{R_{\mathfrak{p}}}(\kappa(\mathfrak{p}),
(H^{d-1}_{I}(M)/\Gamma_{\mathfrak{m}}(H^{d-1}_{I}(M)))_{\mathfrak{p}}),$$is
finite,
where $\kappa(\mathfrak{p})=R_{\mathfrak{p}}/\mathfrak{p}R_{\mathfrak{p}}$.\qed\\

The following lemma and proposition will be sueful in the proof of the next main result of this paper.
An $R$-module $M$ is said to be {\it weakly Laskerian} if the set of associated primes of any quotient
module of $M$ is finite (cf. \cite {DMa, DMa2}).\\

\begin{lem}
Let $R$ be a Noetherian ring and $I$ an ideal of $R$. Then, for any
$R$-module $T$, the following conditions are equivalent:

{\rm(i)} ${\rm Ext}_R^n(R/I,T)$ is weakly Laskerian for all
$n\geq0$.

{\rm(ii)} For any finitely generated $R$-module $N$ with support in
$V(I)$, ${\rm Ext}_R^n(N,T)$ is weakly Laskerian for all $n\geq0$.
\end{lem}
\proof For proving the assertion we use the proof of \cite[Lemma 1]{Ka} and \cite [Lemma 2.2]{DMa2}.\qed\\

\begin{prop}
Let $(R,\mathfrak{m})$ be a local (Noetherian) ring and $M$ a
finitely generated $R$-module. Let $I\subseteq J$ be two ideals of
$R$ such that $\dim R/I=2$. Then, for all $i,j\geq0$, there exists a
finitely generated submodule $L$ of $\Ext^j_R(R/J,H^i_I(M))$ such
that $\Supp (\Ext^j_R(R/J,H^i_I(M))/L)\subseteq V({\mathfrak m}).$
\end{prop}
\proof In view of \cite [Corollary 3.3]{BN3} and Lemma 2.10 the set
$\Ass_R\Ext^j_R(R/J,H^i_I(M))$ is finite. Let
$$\Ass_{R}\Ext^j_R(R/J,H^i_I(M))\backslash
\{\mathfrak{m}\}=\{\mathfrak{p}_{1},\dots,\mathfrak{p}_{n}\}.$$
Since $\mathfrak{p}_{k}\neq \mathfrak{m}$ for all $k=1,\dots,n$ and
$\height\mathfrak{m}/I=2$, it follows that
$\height\mathfrak{p}_k/I\leq1$, and so $\dim
R_{\mathfrak{p}_{k}}/IR_{\mathfrak{p}_{k}}\leq1$ for all
$k=1,\dots,n$. Hence by \cite [Corollary 2.7]{BN3} the
$R_{\mathfrak{p}_{k}}$-modules of
$(\Ext^j_R(R/I,H^i_I(M)))_{\mathfrak{p}_{k}}$ are finitely generated
for all $k=1,\dots,n$. Therefore, by \cite [Lemma 1]{Ka}, the
$R_{\mathfrak{p}_{k}}$-modules of
$(\Ext^j_R(R/J,H^i_I(M)))_{\mathfrak{p}_{k}}$ are finitely generated
for all $k=1,\dots,n$. Consequently for every $k=1,\dots,n$, there
exists a finitely generated submodule $T_{k}$ of
$\Ext^j_R(R/J,H^i_I(M))$ such that $(T_{k})_{\frak p_{k}}=
(\Ext^j_R(R/J,H^i_I(M)))_{\frak p_{k}}.$ Let $K:= T_{1}+\cdots+
T_{n}$. Then $K$ is a finitely generated submodule of
$\Ext^j_R(R/J,H^i_I(M))$ and
$$\Supp (\Ext^j_R(R/J,H^i_I(M))/K) \cap \{\mathfrak{p}_{1},\dots,\mathfrak{p}_{n}\}=\emptyset.$$
 So that $\dim  \Ext^j_R(R/J,H^i_I(M))/K\leq 1$. Now if $\Supp ( \Ext^j_R(R/J,H^i_I(M))/K)\subseteq
 \{\mathfrak{m}\}$, then the result follows for $L=K$. We may
 therefore assume that $$\dim \Ext^j_R(R/J,H^i_I(M))/K=1.$$Let $\Ass_{R}\Ext^j_R(R/J,H^i_I(M))/K
 \backslash \{\mathfrak{m}\}=\{\mathfrak{q}_{1},\dots,\mathfrak{q}_{s}\}$. Then, as
 the above, the $R_{\mathfrak{q}_{t}}$-module of $(\Ext^j_R(R/J,H^i_I(M))/K)_{\mathfrak{q}_{t}}$
 is finitely generated for $t=1,\dots,s$. Hence there exists a finitely generated submodule $L_t/K$ of
 $\Ext^j_R(R/J,H^i_I(M))/K$ such that
$$(\Ext^j_R(R/J,H^i_I(M))/K)_{\frak q_{t}}= (L_t)_{\frak q_{t}}.$$
Let $L:=L_1+\cdots+L_s+K$. Then $L$ is a finitely generated
submodule of $\Ext^j_R(R/J,H^i_I(M))$ and
$$\Supp (\Ext^j_R(R/J,H^i_I(M))/L) \cap \{\mathfrak{q}_{1},...,\mathfrak{q}_{s}\}=\emptyset.$$
Since $\Supp \Ext^j_R(R/J,H^i_I(M))/L\subseteq \Supp \Ext^j_R(R/J,H^i_I(M))/K$,
it follows that $$\Supp \Ext^j_R(R/J,H^i_I(M))/L\subseteq \{\mathfrak{m}\},$$ as required. \qed\\

We are now ready to state and prove the second main theorem of the paper,
which is a characterization of the finiteness of the Bass numbers of ${(d-1)}^{th}$
local cohomology of $M$ with respect to an arbitrary ideal $I$ in terms of the
$I$-cominimaxness of $H^{d-1}_{I}(M)$, where $d=\dim M$. Before we state theorem, let
us recall that a sequence $x_1, \dots, x_n$ of element of an ideal $I$ of $R$ is called
an $I$-filter regular $M$-sequence, provided:
$$\Supp((x_1, \dots, x_{i-1})M:_M x_i )/(x_1, \dots, x_{i-1})M)\subseteq
 V(I),$$
for all $i=1,\dots, n$. In the case of $(R,\frak m)$  a local ring and $I=\frak m$, we recover filter regular
sequences as studied in several papers (cf. e.g. \cite {STC}).

\begin{thm}
Let $(R,\mathfrak{m})$ be a local (Noetherian) ring, $I$ an ideal of
$R$ and $M$ a finitely generated $R$-module of dimension $d$. Then
the following conditions are equivalent:

{\rm(i)} The Bass numbers of the $R$-module $H^{d-1}_{I}(M)$ are
finite.

{\rm(ii)} ${\rm Soc}(H^{d-1}_{I}(M))$ is finitely generated.

{\rm(iii)} $H^{d-1}_{I}(M)$ is $I$-cominimax.
\end{thm}

\proof The conclusion  (i)$\Longrightarrow$(ii) is obviously true.
In order to prove (ii)$\Longrightarrow$ (i), let ${\rm
Soc}\,H^{d-1}_{I}(M)$ is a finitely generated $R$-module.
 Then, according to Melkersson's Theorem \cite [Proposition 1.4]{Me3}, $\Gamma_{\mathfrak{m}}(H^{d-1}_{I}(M))$ is
 Artinian. Therefore the Bass numbers of $\Gamma_{\mathfrak{m}}(H^{d-1}_{I}(M))$ are finite.
 Now, it follows from Theorem 2.9 that the Bass numbers of $H^{d-1}_{I}(M)$ are finite, as required.\\
In order to prove the implication (iii)$\Longrightarrow$(ii), it
follows from \cite [Corollary 2.8]{ANV} that the $R$-module
$\Ext^{j}_{R}(R/\mathfrak{m},H^{d-1}_{I}(M))$ is minimax for all
$j\geq0$. Now it easily seen from definition that
$\Ext^{j}_{R}(R/\mathfrak{m},H^{d-1}_{I}(M))$ is of finite length
for all $j\geq0$. In particular the $R$-module
$\Hom_{R}(R/\mathfrak{m},H^{d-1}_{I}(M))$
 is finitely generated, as required.

Finally for the proof of (i)$\Longrightarrow$ (iii), let $\dim
R/I=n$. For $n\leq 1$ this follows by \cite [Corollary 2.7]{BN3}.
Hence we may assume that $n\geq 2$. Let $x_{1},...,x_{n-1}\in
\mathfrak{m}$ be a filter regular sequence on $R/I$. Let $K_{s}:=I+(x_{1},\dots,x_{s})$ for
$s=0,1,\dots,n-1$. We first show by induction on $t$, that the
$R$-module $\Ext^{j}_{R}(R/K_{n-t-1},H^{d-1}_{I}(M))$ is minimax for
all $j\geq0$. Since $x_{1}+I,\dots,x_{n-1}+I$ is a part of system of
parameters for $R/I$, it follows that $\dim R/K_{n-1}=1 $, and so
the result is evidently true for $t=0$, by Theorem 2.3. Suppose that
$n-1\geq t>0$ and the case $t-1$ is settled. Since
$\Gamma_{\mathfrak{m}}(R/K_{n-t-1})$ is a finitely generated
$R$-module with $\Supp (\Gamma_{\mathfrak{m}}(R/K_{n-t-1}))\subseteq
\{\mathfrak{m}\}$, it follows that from \cite [Lemma 1]{Ka} and
hypothesis (i) that the $R$-module
$\Ext^{j}_{R}(\Gamma_{\mathfrak{m}}(R/K_{n-t-1}),H^{d-1}_{I}(M))$ is
finitely generated for all $j\geq0$. Now, as
$\Gamma_{\mathfrak{m}}(R/K_{n-t-1})=J/K_{n-t-1}$ for some ideal $J$
of $R$, from the exact sequence $$0\longrightarrow J/K_{n-t-1}
\longrightarrow R/K_{n-t-1} \longrightarrow R/J \longrightarrow0,$$
we get the following long exact sequence,
$$\cdots \longrightarrow \Ext^{j-1}_{R}(J/K_{n-t-1},H^{d-1}_{I}(M)) \longrightarrow
\Ext^{j}_{R}(R/J,H^{d-1}_{I}(M)) \longrightarrow$$ $$\longrightarrow
\Ext^{j}_{R}(R/K_{n-t-1},H^{d-1}_{I}(M)) \longrightarrow
\Ext^{j}_{R}(J/K_{n-t-1},H^{d-1}_{I}(M)) \longrightarrow \cdots.$$
Hence we get that the $R$-module
$\Ext^{j}_{R}(R/K_{n-t-1},H^{d-1}_{I}(M))$ is minimax if and only if
the $R$-module
 $\Ext^{j}_{R}(R/J,H^{d-1}_{I}(M))$ is minimax for all $j\geq 0$. Now, as $x_{1},\dots,x_{n-1}$ is a
 filter regular sequence on $R/I$ in $\mathfrak{m}$, it follows from definition that
$$x_{n-t}\not\in \bigcup_{\mathfrak{p}\in \Ass_{R}(R/K_{n-t-1})\backslash
\{\mathfrak{m}\}}\mathfrak{p}.$$ Therefore, it follows from
$$\Ass_{R}(R/J)=\Ass_{R}(R/K_{n-t-1})\backslash \{\mathfrak{m}\}$$
 that $x_{n-t}\not\in \cup_{\mathfrak{p}\in
\Ass_{R}R/J}\mathfrak{p}$, and so $x_{n-t}$ is an $R/J$-regular
element in $\mathfrak{m}$. Consequently, the exactness
$$0\longrightarrow R/J\stackrel{x_{n-t}}\longrightarrow R/J\longrightarrow
R/J+Rx_{n-t} \longrightarrow 0,$$ implies that the sequence
$$0\longrightarrow
\Ext^{j-1}_{R}(R/J,H^{d-1}_{I}(M))/x_{n-t}\Ext^{j-1}_{R}(R/J,H^{d-1}_{I}(M))
\longrightarrow
$$$$\longrightarrow \Ext^{j}_{R}(R/J+Rx_{n-t},H^{d-1}_{I}(M))
\longrightarrow (0:_{\Ext^{j}_{R}(R/J,H^{d-1}_{I}(M))}x_{n-t})
\longrightarrow 0,$$ is exact. Since $\Supp (R/J+Rx_{n-t})\subseteq
V(K_{n-t})$, (note that $K_{n-t}\subseteq J+Rx_{n-t}$,) the
induction hypothesis and \cite [Lemma 1]{Ka} yield that the
$R$-module $$\Ext^{j}_{R}(R/J+Rx_{n-t},H^{d-1}_{I}(M)),$$ is minimax,
for all $j\geq0$. Therefore the $R$-module
$(0:_{\Ext^{j}_{R}(R/J,H^{d-1}_{I}(M))}x_{n-t})$ is minimax. Now,
let $\mathfrak{p}\in \Supp
H^{d-1}_{I}(M)\backslash\{\mathfrak{m}\}$. Then, by \cite
[Proposition 5.1]{Me3} and Vanishing theorem the
$R_{\mathfrak{p}}$-module $(H^{d-1}_{I}(M))_{\mathfrak{p}}$, is
$IR_{\mathfrak{p}}$-cofinite. As $I\subseteq J \subseteq
J+Rx_{n-t}$, it follows from \cite [Lemma 1]{Ka} that the
$R_{\mathfrak{p}}$-module
$(\Ext^{j}_{R}(R/J,H^{d-1}_{I}(M)))_{\mathfrak{p}}$ is finitely
generated. Moreover, as in view of \cite [Corollary 3.3]{NS} the set
$\Supp H^{d-1}_I(M)$ is finite, it follows that the set
$\Ass_R\Ext^j_R(R/J,H^{d-1}_I(M))$ is finite. Therefore, by the
proof of Lemma 2.11, there exists a finitely generated submodule $L$
of $\Ext^j_R(R/J,H^{d-1}_I(M))$ such that
$$\Supp\Ext^j_R(R/J,H^{d-1}_I(M))/L\subseteq \{\mathfrak{m}\}.$$
Next, from the exact sequence
$$0\longrightarrow L \longrightarrow \Ext^j_R(R/J,H^{d-1}_I(M))\longrightarrow \Ext^j_R(R/J,H^{d-1}_I(M))/L \longrightarrow0,$$
we get the exact sequence $$(0:_{\Ext^j_R(R/J,H^{d-1}_I(M))}x_{n-t})
\longrightarrow (0:_{\Ext^j_R(R/J,H^{d-1}_I(M))/L}x_{n-t})
\longrightarrow Ext^{1}_{R}(R/Rx_{n-t},L).$$ Hence, it yields that
the $R$-module $(0:_{\Ext^j_R(R/J,H^{d-1}_I(M))/L}x_{n-t})$ is
minimax. Since $\Supp
(0:_{\Ext^j_R(R/J,H^{d-1}_I(M))/L}x_{n-t})\subseteq
\{\mathfrak{m}\}$, it follows from the definition that the
$R$-module $(0:_{\Ext^j_R(R/J,H^{d-1}_I(M))/L}x_{n-t})$ is Artinian.
As $(0:_{\Ext^j_R(R/J,H^{d-1}_I(M))/L}x_{n-t})$ is
$Rx_{n-t}-$torsion,  according to Melkersson \cite
[Proposition 1.4]{Me3} $\Ext^j_R(R/J,H^{d-1}_I(M))/L$ is an Artinian
$R$-module. That is $\Ext^j_R(R/J,H^{d-1}_I(M))$ is a minimax
$R$-module, as
required.\qed\\

\begin{cor}
Let $(R,\mathfrak{m})$ be a regular local ring of dimension $d\geq
1$ containing a field and let $I$ be an ideal of $R$. Then the
$R$-module $H^{d-1}_{I}(R)$ is $I$-cominimax .
\end{cor}

\proof The result follows from Theorem 2.12 and the fact that the
Bass numbers of the $R$-module $H^{d-1}_{I}(R)$ are finite (see
\cite {HS, Ly2}).\qed\\

\begin{cor}
Let $(R,\mathfrak{m})$ be a local (Noetherian) ring and $I$ an ideal
of $R$. Let $M$ be a finitely generated $R$-module of dimension $d$.
Then the Bass numbers of the $R$-module $H^{d-1}_{I}(M)$ are finite,
whenever $\mathfrak{m}\not\in \Ass_{R}H^{d-1}_{I}(M)$.
\end{cor}
\proof Since $\Ass_R({\rm Soc}\, H^{d-1}_{I}(M))=
V(\mathfrak{m})\cap\Ass_R H^{d-1}_{I}(M)$, it follows from
$$\mathfrak{m}\not\in \Ass_{R}(H^{d-1}_{I}(M))$$ that $\Ass_R({\rm
Soc}\, H^{d-1}_{I}(M))=\emptyset$. Hence ${\rm Soc}\,
H^{d-1}_{I}(M)=0$,
and so the result now follows from Theorem 2.12.\qed\\

\begin{cor}
Let $(R,\frak{m})$ be a  local (Noetherian) ring. Let  $I$ be
an ideal of $R$ and $x\in {\frak m}$ such that  $I\subseteq Rx$.
Then for each non-zero finitely generated $R$-module $M$ of dimension
$d\geq 1$, the $R$-module $H^{d-1}_{I}(M)$ is $I$-cominimax.
In particular, the Bass numbers of $H^{d-1}_{I}(M)$ are finite.
\end{cor}

\proof When $d=1$ there is nothing to prove. In the case that $d=2$,
it follows from \cite[Theorem 2.3]{BN1} that the $R$-module
${\rm Hom}_R(R/I, H^{1}_I(M))$ is finitely
generated and so the $R$-module ${\rm Soc}(H^{1}_I(M))$ is finitely
generated. Now the assertion follows from Theorem 2.12. Finally, for
$d\geq 3$, it follows from \cite[Lemma 2.5]{AB} that the map
$H^{d-1}_I(M)\stackrel{x} \longrightarrow H^{d-1}_I(M)$ is an
isomorphism, and so ${\rm Soc}(H^{d-1}_I(M))=0$. Now the assertion
again follows from Theorem 2.12.\qed\\

\begin{cor}
Let $(R,\mathfrak{m})$ be a regular local ring and $I$ an ideal of $R$ such that ${\rm height}(I)=1$. Then for
each non-zero finitely generated $R$-module $M$ of dimension $d\geq
1$, the $R$-module $H^{d-1}_{I}(M)$ is $I$-cominimax. In particular,
the Bass numbers of $H^{d-1}_{I}(M)$ are finite.
\end{cor}

\proof Since $R$ is a UFD and ${\rm height}(I)=1$, it follows that
$I$ is contained in a principal prime ideal. Now
the assertion follows from Corollary 2.15.\qed\\

We are now ready to show that the subjects of the previous results can be
used to prove a  characterization of the finiteness of the Bass numbers of $i^{th}$
local cohomology of $M$ with respect to an  ideal $I$  of dimension 2, i.e., where $\dim R/I=2$,
in terms of the $I$-cominimaxness of $H^{i}_{I}(M)$. The main result is Theorem 2.19.
The following theorem will serve to shorten the proof of the  main theorem.\\

\begin{thm}
Let $(R,\mathfrak{m})$ be a local (Noetherian ) ring and let $I$ be
an ideal of $R$ with $\dim R/I=2$. Let $M$ be a finitely generated
$R$-module and $i$ a non-negative integer. Then the following
conditions are equivalent:

{\rm(i)} The Bass numbers of $H^{i}_{I}(M)$ are finite.

{\rm(ii)} The $R$-module $H^{i}_{I}(M)$ is $I$-cominimax.
\end{thm}

\proof First we show the implication (i)$\Longrightarrow$ (ii).
Suppose that $\Gamma_{\mathfrak{m}}(R/I):=J/I$. Then $\Supp
J/I\subseteq \{\mathfrak{m}\}$, and so it follows from \cite [Lemma
1]{Ka} and the assumption (i) that the $R$-module
$\Ext^{j}_{R}(J/I,H^{i}_{I}(M))$ is finitely generated for all
$j\geq 0$. Now, the exact sequence $$0\longrightarrow J/I
\longrightarrow R/I \longrightarrow R/J \longrightarrow0$$ induces a
long exact sequence
$$\cdots \longrightarrow \Ext^{j-1}_{R}(J/I,H^{i}_{I}(M)) \longrightarrow
\Ext^{j}_{R}(R/J,H^{i}_{I}(M)) \longrightarrow$$$$\longrightarrow
\Ext^{j}_{R}(R/I,H^{i}_{I}(M)) \longrightarrow
\Ext^{j}_{R}(J/I,H^{i}_{I}(M)) \longrightarrow \cdots.$$ Therefore
$H^{i}_{I}(M)$ is $I$-cominimax if and only if the $R$-module
$\Ext^j_R(R/J,H^{i}_I(M))$ is minimax for all $j\geq 0$.

It is easy to see that $\Gamma_{\mathfrak{m}}(R/J)=0$. Hence, by
\cite [Lemma 2.1.1]{BS} $\mathfrak{m}$ contains an $R/J$-regular
element $x$. Then $\dim R/J/x(R/J)=\dim R/J-1$, and so $\dim R/J+Rx=
1$. Furthermore, the exact sequence $$0\longrightarrow
R/J\stackrel{x}\longrightarrow R/J\longrightarrow R/J+Rx
\longrightarrow0$$ gives us an exact sequence
$$0\longrightarrow
\Ext^{j-1}_{R}(R/J,H^{i}_{I}(M))/x\Ext^{j-1}_{R}(R/J,H^{i}_{I}(M))
\longrightarrow $$$$\longrightarrow
\Ext^{j}_{R}(R/J+Rx,H^{i}_{I}(M)) \longrightarrow
(0:_{\Ext^{j}_{R}(R/J,H^{i}_{I}(M))}x) \longrightarrow0$$  Now,
since $\dim R/J+Rx=1$, it follows from the assumption (i) and
Theorem 2.3 that the $R$-module $\Ext^{j}_{R}(R/J+Rx,H^{i}_{I}(M))$
is minimax, for all $j\geq0$. On the other hand, in view of
Proposition 2.11, there exists a finitely generated submodule $L$ of
$\Ext^{j}_{R}(R/J,H^{i}_{I}(M))$ such that $\Supp
(\Ext^{j}_{R}(R/J,H^{i}_{I}(M))/L)\subseteq \{\mathfrak{m}\}$. Now,
from the exact sequence$$0\longrightarrow L \longrightarrow
\Ext^j_R(R/J, H^{i}_I(M))\longrightarrow \Ext^j_R(R/J, H^{i}_I(M))/L
\longrightarrow0$$ we obtain the exact sequence
$$(0:_{\Ext^j_R(R/J, H^{i}_I(M))}x) \longrightarrow
(0:_{\Ext^j_R(R/J, H^{i}_I(M))/L}x) \longrightarrow
\Ext^{1}_{R}(R/Rx, L),$$ which implies that the $R$-module
$(0:_{\Ext^j_R(R/J, H^{i}_I(M))/L}x)$ is minimax. As $$\Supp
(0:_{\Ext^j_R(R/J, H^{i}_I(M))/L}x)\subseteq \{\mathfrak{m}\},$$ it
follows from the definition that the $R$-module
$(0:_{\Ext^j_R(R/J,H^{i}_I(M))/L}x)$ is Artinian. Since
$(0:_{\Ext^j_R(R/J,H^{i}_I(M))/L}x)$ is $Rx$-torsion, it follows
from Melkersson's result (see \cite [Proposition 1.4]{Me3}) that
$\Ext^j_R(R/J,H^{i}_I(M))/L$ is an Artinian $R$-module, and so
$\Ext^j_R(R/J,H^{i}_I(M))$ is a minimax $R$-module, as
required.\\
For prove (ii)$\Longrightarrow$(i) first observe that the $R$-module
$\Ext^{j}_{R}(R/\mathfrak{m},H^{i}_{I}(M))$ is minimax, for all
$j\geq0$, and so by the definition,
$\Ext^{j}_{R}(R/\mathfrak{m},H^{i}_{I}(M))$ has finite length.
Therefore $\mu^{j}(\mathfrak{m},H^{i}_{I}(M))$ is finite, for all
$j\geq0$. Now, let $\mathfrak{p}\in \Supp H^{i}_{I}(M)\backslash
V(\mathfrak{m})$. Then $\dim R_{\mathfrak{p}}/IR_{\mathfrak{p}}\leq
1$ as easily seen. Therefore the $R_{\mathfrak{p}}$-module of
$H^{i}_{I}(M)_{\mathfrak{p}}$ is $IR_{\mathfrak{p}}$-cofinite, and
so the Bass numbers of $H^{i}_{I}(M)_{\mathfrak{p}}$ are finite.
Since the Bass numbers are stable under the localization, it follows
that the Bass numbers of $H^{i}_{I}(M)$ are finite, as
required.\qed\\

The next result is an immediate consequence of Theorem 2.17 and Huneke-Sharp-Lyubeznik' results.\\

\begin{cor}
Let $(R,\mathfrak{m})$ be a regular local ring of dimension $d\geq
2$ containing a field and let $I$ be an ideal of $R$ with $\dim
R/I=2$.
 Then the $R$-module $H^{i}_{I}(R)$ is $I$-cominimax for all $i\geq0$.
\end{cor}

\proof The result follows from Theorem 2.17 and the fact that the
Bass numbers of the $R$-module $H^{i}_{I}(R)$ are finite  (see \cite
{HS, Ly2}), for all $i\geq0$.\qed\\

We are now ready to state and prove the third main theorem of the  paper,
which is a characterization of the finiteness of the Bass numbers of $i^{th}$
local cohomology of $M$ with respect to an  ideal $I$  of dimension 2,
in terms of the $I$-cominimaxness of $H^{i}_{I}(M)$.\\

\begin{thm}
Let $(R,\mathfrak{m})$ be a local (Noetherian) ring, and let $I$ be
an ideal of $R$ such that $\dim R/I=2$. Let $M$ be a finitely
generated $R$-module and $t$ a non-negative integer. Then the
following conditions are equivalent:

{\rm(i)} ${\rm Soc} H^{i}_{I}(M)$ is finitely generated for all
$i\leq t$.

{\rm(ii)} $\Ext^{j}_{R}(R/\mathfrak{m},H^{i}_{I}(M))$ is finitely
generated for all $j \geq 0$ and for all $i\leq t-1$.

{\rm(iii)} $H^{i}_{I}(M)$ is $I$-cominimax, for all $i\leq t-1$.
\end{thm}

\proof The implication (ii)$\Longrightarrow$ (i) follows from
\cite[Corollary 3.5]{Kh}. In order to prove (i)$\Longrightarrow$
(ii), as $\Ass_R H^{i}_{I}(M)$ is finite (see \cite[Corollary
4.3]{BN3}), we can find an element $x\in \mathfrak{m}$, such that $x
\not \in \bigcup_{\mathfrak{p}\in {\rm Assh}(R/I)}\mathfrak{p},$ and
$V(Rx) \cap
(\bigcup_{i=0}^{t}\Ass_{R}H^{i}_{I}(M))\subseteq V(\mathfrak{m}).$\\
Now, let $L:= I+Rx$. Then it is easy to see that $\dim R/L=1$, and
$$H^{0}_{L}(H^{i}_{I}(M))=H^{0}_{Rx}(H^{i}_{I}(M))=H^{0}_{\mathfrak{m}}(H^{i}_{I}(M)),$$
is an Artinian $R$-module for all $i\leq t$, by \cite[Proposition
4.1]{Me3} and assumption (i).  Now, for all $i\leq t-1$, there
exists the short exact sequence
$$0\longrightarrow H^{1}_{Rx}(H^{i}_{I}(M)) \longrightarrow H^{i+1}_{L}(M) \longrightarrow
H^{0}_{Rx}(H^{i+1}_{I}(M))\longrightarrow 0,$$(cf. for instance
\cite[Corollary 3.5]{Sc}).
 Since $H^{i+1}_{L}(M)$ is $L$-cofinite and
 $$H^{0}_{Rx}(H^{i+1}_{I}(M))=H^{0}_{\mathfrak{m}}(H^{i+1}_{I}(M))$$
  is Artinian, it follows from \cite[Lemma 2.1]{BN3}
  that $\Ext^{j}_{R}(R/L,H^{1}_{Rx}(H^{i}_{I}(M)))$ is minimax for
  all $j\geq0$ and $i\leq t-1$. Furthermore, since $H^{k}_{Rx}(H^{i}_{I}(M))=H^{k}_{L}(H^{i}_{I}(M))$,
  for all $k\geq0$, it follows that $H^{k}_{L}(H^{i}_{I}(M))=0$ for all $k\geq 2$.
  Hence $\Ext^{j}_{R}(R/L,H^{k}_{L}(H^{i}_{I}(M)))$ is minimax for all $j\geq 0$, $k\geq 0$ and $i\leq t-1$,
  and so in view of \cite[Proposition 3.9]{Me3}, the $R$-module $\Ext^{j}_{R}(R/L,H^{i}_{I}(M))$ is minimax
  for all $j \geq 0$ and $i\leq t-1$. Consequently, by \cite[Theorem 2.7]{ANV}
   the $R$-module $K_{i,\,j}:=\Ext^{j}_{R}(R/\mathfrak{m},H^{i}_{I}(M))$ is minimax
   for all $j \geq 0$ and $i\leq t-1$. Hence there is a finitely generated submodule
   $N_{i,\,j}$ of $K_{i,\,j}$, such that $K_{i,j}/N_{i,\,j}$ is Artinian.
   As $\mathfrak{m}K_{i,\,j}=0$ and
   $\mathfrak{m}(K_{i,\,j}/N_{i,\,j})=0$, it follows that $K_{i,\,j}$ has finite length, as required.

 For proving (ii)$\Longrightarrow$(iii), first we use the proof of implication (ii)$\Longrightarrow$(i) in Theorem
 2.17, to see that the Bass numbers of $H^{i}_{I}(M)$ are finite, for all $i\leq t-1$.
Now, the assertion follows from Theorem 2.17. Finally the
implication (iii)$\Longrightarrow$(ii) follows again from Theorem 2.17.\qed\\

We end the paper by drawing some concequences of Theorem 2.19. The first
consequence will help us to construct a certain subset of $\Spec R$ which is countable.\\

\begin{cor}
Let $(R,\frak{m})$ be a local (Noetherian) ring and let $I$ be
an ideal of $R$  with $\dim R/I=2$. Let $M$ be a finitely
generated $R$-module such that  $\frak{m} \not \in \cup_{i\geq
0}\Ass_R(H^i_I(M))$. Then
 $\Ext^{j}_{R}(R/\mathfrak{m},H^{i}_{I}(M))$ is finitely
generated for all $j \geq 0$ and all $i\geq 0$.
\end{cor}

\proof Since $\frak{m}  \not \in \cup_{i\geq 0}\Ass_R(H^i_I(M))$,
it follows that ${\rm Soc} H^{i}_{I}(M)=0$  for all $i\geq 0$.
Now the assertion follows  from Theorem 2.19.\qed\\

\begin{prop}
Let $R$ be a Noetherian ring and $I$ an ideal of $R$ such
that $\dim R/I=2$. For a finitely generated $R$-module $M$, let
$$\Sigma:=\{\frak p \in {\rm Spec}(R)\,\,:\,\,\mu^j(\frak p,H^{i}_{I}(M)
)=\infty\,\,{\rm for}\,\,{\rm some}\,\,{\rm integers}\,\,i\geq
0\,\,{\rm and}\,\,j\geq 0 \}.$$ Then $\Sigma$ is  countable  and
$\Sigma \subseteq {\rm Max}(R)$.
\end{prop}

\proof

In order to show that $\Sigma \subseteq {\rm Max}(R)$,
let $\frak p \in \Sigma\backslash {\rm Max}(R)$. Then it follows from the definition that $\frak
p \in V(I)$, and it is easy to see that ${\rm height}(\frak p/I)\leq 1$. Hence  it
follows from \cite[Corollary 2.10]{BN3} that $\frak p\not \in\Sigma $, which is a
contradiction. Thus $\Sigma\subseteq {\rm Max}(R)$.

On the other hand, as for all integers $n\geq 1$ and $i\geq0$, the
sets ${\rm Ass}_R{\rm Ext}^i_R(R/I^n, M)$ are finite  and $H^i_I(M)
= \underset{n\geq1} {\varinjlim}\,\, \text{Ext}^i_R(R/I^n, M),$ it
follows that the set $\bigcup_{i\geq 0}\Ass_RH^i_I(M)$ is countable.
Now, in order to complete the proof, it is enough to show that
$\Sigma\subseteq \bigcup_{i\geq 0}\Ass_RH^i_I(M).$ To do this
suppose that the contrary is true. Then there is $\frak m \in
\Sigma\backslash \bigcup_{i\geq 0}\Ass_RH^i_I(M)$. Hence, by
\cite[Corollary 2.10]{BN3} we have ${\rm height}(\frak m/I)= 2$, and
so it follows from Corollary 2.20 that
$\frak m \not \in \Sigma$, which is a contradiction.\qed\\

\begin{cor}
Let $R$ be a Noetherian ring, $M$ a finitely generated $R$-module
and $I$  an ideal of $R$ such that $\dim M/IM=2$ . Let
$$\Sigma:=\{\frak p \in {\rm Spec}(R)\,\,:\,\,\mu^j(\frak p,H^{i}_{I}(M)
)=\infty\,\,{\rm for}\,\,{\rm some}\,\,{\rm integers}\,\,i\geq
0\,\,{\rm and}\,\,j\geq 0 \}.$$ Then $\Sigma$ is  countable  and
$\Sigma \subseteq {\rm Max}(R).$
\end{cor}

\proof

As $\dim M/IM=2$ it follows from that
$\dim R/I+{\rm Ann}_RM=2$. Now the result follows from
 $H^{i}_{I}(M)\cong H^{i}_{I+{\rm Ann}_R(M)}(M)$ and Proposition 2.21. \qed\\

\begin{cor}
Let $(R,\mathfrak{m})$ be a local (Noetherian) ring, and let $I$ be
an ideal of $R$ such that $\dim R/I=3$.  Let $M$ be a finitely
generated $R$-module. Then the set
$$\Sigma:=\{\frak p \in {\rm Spec}(R)\,\,:\,\,\mu^j(\frak p,H^{i}_{I}(M)
)=\infty\,\,{\rm for}\,\,{\rm some}\,\,{\rm integers}\,\,i\geq
0\,\,{\rm and}\,\,j\geq 0 \},$$ is  countable.
\end{cor}

\proof Let $\frak p \in \Sigma\backslash \{\frak m\}$. Then by
\cite[Corollary 2.10]{BN3} we have ${\rm height}(\frak p/I)= 2$,
and so it follows from Corollary 2.20 that $\frak p \in
\bigcup_{i\geq 0}\Ass_RH^i_I(M)$. Thus
 $\Sigma\backslash \{\frak m\}\subseteq \bigcup_{i\geq 0}\Ass_RH^i_I(M)$,
and so the set $\Sigma$ is  countable.  \qed\\

\begin{cor}
Let $(R,\mathfrak{m})$ be a local (Noetherian) ring, $M$  a
finitely generated $R$-module and $I$ an ideal of $R$  such that
$\dim M/IM=3$. Then the set
$$\Sigma:=\{\frak p \in {\rm Spec}(R)\,\,:\,\,\mu^j(\frak p,H^{i}_{I}(M)
)=\infty\,\,{\rm for}\,\,{\rm some}\,\,{\rm integers}\,\,i\geq
0\,\,{\rm and}\,\,j\geq 0 \}$$ is  countable.
\end{cor}

\proof As $\dim M/IM=3$ it follows that $\dim R/I+{\rm Ann}_R(M)=3$.
Now the result follows from  $H^{i}_{I}(M)\cong H^{i}_{I+{\rm
Ann}_R(M)}(M)$ and Corollary 2.23
.\qed\\

\begin{center}
{\bf Acknowledgments}
\end{center}
The authors are deeply grateful to the referee for a very careful
reading of the manuscript and many valuable suggestions. We also would like to thank Professors Hossein Zakeri and Kamran
Divaani-Aazar for their careful reading of the first draft and many helpful suggestions.

\end{document}